\documentclass[]{article}
\usepackage{amsmath}
\usepackage{amssymb}
\usepackage{color}
\usepackage{graphicx}
\usepackage{booktabs}
\newcommand{\RM}[1]{\MakeUppercase{\romannumeral #1}}

\newcommand{\WW}{W} 
\newcommand{\no}{\nonumber} 
\begin{document}
\title{Finite-difference-based simulation and adjoint optimization of gas networks\protect}

\author{Sonja Hossbach\footnote{Technische Universität Berlin, Numerische Fluiddynamik,  Mueller-Breslau-Str.~15,
		10623 Berlin,  {e-mail: sonja.hossbach@tnt.tu-berlin.de}} , Mathias Lemke, Julius Reiss}

\maketitle

\abstract{
The stable operation of gas networks is an important optimization target. 
While for this task commonly finite volume methods are used, we introduce a new finite difference approach. 
With a summation by part formulation for the spatial discretization, we get well-defined fluxes between the pipes.
This allows a simple and explicit formulation of the coupling conditions at the node.
From that, we derive the adjoint equations for the network simply and transparently. 
The resulting direct and adjoint equations are numerically efficient and easy to implement.
The approach is demonstrated by the optimization of two sample gas networks.}





\section{Introduction}

Gas networks play a crucial role in the  energy supply.
By this, their stability  is of great importance.
However, controlling gas networks with a variety of sources is a major challenge due to unpredictable time-dependent demands from end-users.
This is exacerbated by the advent of renewable energy sources such as wind power and solar systems as gas power plants are used to compensate fluctuations from other energy sources and power-to-gas is an option for storing excess renewable energy.
Thus, the control of gas transport is becoming more and more complex.
To achieve stable operation, gas networks must be modeled and simulated in terms of a time-dependent transient technical optimization.
The main objective is to guarantee the security of supply and the system stability of the gas networks.

Early gas network simulations used steady-state methods \cite{osiadacz1988method}. With the increasing complexity of the networks for most applications an unsteady simulation is needed today.
Nevertheless, the steady-state analysis is still used in some cases \cite{bermudez2015simulation, schmidt2015high}.
Complex networks with compressor stations and valves need operation decisions which leads to mixed-integer optimization problems.
A comparison of different approaches (NLP - nonlinear programming, MILP - mixed-integer linear programming, MINLP - mixed-integer nonlinear programs) is given in Pfetsch et al. \cite{pfetsch2015validation}.

Also, the coupling conditions are the aim of current research. 
Equal pressure, equal dynamic pressure, or equal momentum are discussed
\cite{herty2008coupling}.
Domschke et al.\cite{DomschkeHillerLangTischendorf2017} and Lang et al.\cite{LangLeugeringMartinTischendorf2017} give an introduction to the topic of modeling and optimizing gas networks and show the different levels of complexity that need to be considered.
One key challenge is the expensive solution of Riemann problems in the typically applied finite volume (FV) method.

The adjoint approach is a powerful and useful tool of functional analysis.
In the field of fluid mechanics, it has been used for stability analysis and control as well as optimization purposes \cite{Jameson1988,GilesPierce2000}.
The adjoint information is commonly used as a gradient to optimize flow configurations through geometry adjustments \cite{Jameson1995} or active flow control \cite{CarnariusThieleOzkayaNemiliGauger2013,LemkeCitroGiannetti2021}. 
Besides, they are used for data assimilation tasks \cite{YangRobinsonHeitzMemin2015, LemkeSesterhenn2016} and for analyzing and optimizing reactive flow configurations \cite{LemkeReissSesterhenn2014, GrayLemkeReissPaschereitSesterhennMoeck2017}. 
Furthermore, the adjoint approach is used in the field of aeroacoustics \cite{Freund2011, SteinStraubeSesterhennWeinzierlLemke2019}.

Also in the context of gas networks, the adjoint approach is already established. 
It has been used to make decisions in hierarchical models \cite{domschke2011adjoint}, error estimations \cite{domschke2015adjoint} and for optimizing gas networks \cite{o2018optimizing,KolbOliver2011} using a finite volume discretization.
\medskip 

Here, we present a new finite differences (FD) approach to simulate and optimize gas networks in contrast to the commonly used finite volume
methods. The pipes are discretized with a summation by parts (SBP) scheme
guaranteeing well-defined fluxes. This allows implementing a
fully explicit predictor-corrector-method for the coupling in the
junctions. No Riemann problem needs to be solved. 
It yields an efficient and transparent numerical scheme for
the full network.
Furthermore, the adjoint of the network can be easily
constructed from the adjoint of single pipes.

The paper is organized as follows: At first, we propose  the network model with a summation by part finite difference scheme and derive the coupling conditions in section 2. Then, we introduce the adjoint approach for the isothermal Euler equations and extend it to the network formulation in section 3. Both, the finite difference scheme and the adjoint are applied to two networks to show the capability of the approach in section 4.


\section{The Network Model}

\subsection{Physical Description}

Gas networks consist of several different elements such as pipes, nodes, compressors, and valves. 
In this paper, we focus on the introduction of a novel discretization and consider pipes and nodes only. 
The fluid dynamical processes are modeled as isothermal compressible flow. 
This is motivated by the fact that long pipelines with non-perfect isolation are assumed to adapt to the temperature of the surrounding.
The coupling of the pipes is constructed from two constraints at each junction. 
First, the assumption of equal pressure of all pipes at each junction (C1), and second, the conservation of mass (C2).

{We refer to individual  pipes with the superscript $\alpha$ and the nodes or junctions with 
$k = \{1, \dots , K \}$ and $K$ being the total number of nodes in the network.}

\paragraph{Pipes}
The flow in the pipes is modeled by the one-dimensional, isothermal Euler equations: 
\begin{align}
    \partial_t \left(\varrho A\right)^\alpha +       \partial_x    \left(\varrho u A \right)^\alpha  &=0 \label{eq_euler_conti}  \\ 
    \partial_t  \left(\varrho u A \right)^\alpha +       \partial_x \left(\varrho u u A \right)^\alpha + A^\alpha \partial_x p^\alpha  &=0. \label{eq_euler_momentum} 
\end{align}
Therein, $\varrho(x,t)$ denotes the density, 
$u(x,t)$ the velocity and $A$ the the cross-sectional area.
In the following we assume constant cross-section for each pipe, possibly different for each  pipe, with index $\alpha$. 
Due to the isothermal assumption, the pressure in the momentum equations is replaced by $ p = c^2 \varrho $, with the speed of sound $c$, which depends on the temperature and is assumed to be prescribed.
We do not explicitly emphasize that the governing variables are functions of space and time for the sake of brevity.

\paragraph{Junctions}
A junction $k$ is given by prescribing $N_k$, being the set of participating pipes $\alpha$ and the index marking their start ($j=1$) or end ($j=N^\alpha$) point respectively.
To connect several pipes two conditions (C1) and (C2) are applied. 
The first condition (C1) is an equal pressure for all pipes at one junction or equivalent the equal density $\bar \varrho_k $: 
\begin{equation}
    \varrho^\alpha_j\left.\right|_{(\alpha,j )\in N_k} = \bar \varrho_k   .
    \label{equalPressure}
\end{equation}
The second condition (C2) is the conservation of mass at the junction $k$.
\begin{equation}
    \sum \limits_{\substack{(\alpha,j) \in N_k \\ (in)}} \left( \varrho u A\right)^\alpha_j = \sum \limits_{\substack{(\alpha,j) \in N_k \\ (out)}} \left( \varrho u A\right)^\alpha_j  \label{eq:conservation_mass}
\end{equation}
Therein, (in)  is encoding the start ($j=1$) and (out) the end of the pipe ($j=N^\alpha$). 
Thus, equation \eqref{eq:conservation_mass} can be written more compact by introducing 
\begin{equation}
    n_j = 
    \left\{
    \begin{array}{rl}
            -1 & \quad \text{for} \quad j =1   \\
            1  & \quad \text{for} \quad j = N^\alpha   
    \end{array}
    \right.  
\end{equation}
with 
\begin{equation}
  \sum \limits_{(\alpha,j) \in N_k} \left( \varrho u A \right)^\alpha_j n_j=0 .
  \label{Kirchhoff} 
\end{equation}
 $n_j$ can be interpreted as the normal vector pointing outward the pipe.
The expression is symmetric, up to the choice of the direction $x$ and the velocity. 
Changing both keeps the expression identical so that we do not prescribe the flow direction in the following.

While the conservation of mass is evident, the equal pressure condition is a common model, approximating a real junction, used  frequently  \cite{o2018optimizing, KolbOliver2011, BandaHertyKlar2006}.

\subsection{Basics of the discretization - SBP} \label{sec:numerical_implementation}
The use of summation by parts (SBP) derivative matrices is essential for the proposed method. 
These SBP matrices can be used to realize well-defined fluxes as summarized in the following \cite{CarpenterGottliebAbarbanel1994}.
This is the key to define a simple numerical procedure for the coupling in the nodes.  
Different SBP-implementations are available, see Sv\"ard and Nordstr\"om \cite{svard2014review} for a detailed review.
We rely on explicit SBP matrices,  for which various orders are given by Strand \cite{Strand1994}. 
In general, SBP schemes consist of a special spatial difference scheme, so that the discrete derivative  operator  can be decomposed as 
\begin{equation}
D_x = W^{-1} S,   \label{SBPdecomp}
\end{equation}
with $D_x,W,S \in \mathbb{R}^{N,N}$.
The matrix $S$ is skew symmetric beside the elements $S_{1,1}$  and $S_{N,N}$. The transposed of $S$ is 
\begin{equation}
S^T = -S + B \label{SBPtrans}
\end{equation}
with 
$B_{i,j} = -\delta_{i,1}\delta_{j,1} + \delta_{i,N}\delta_{j,N}    $. 
The matrix $W$ is a symmetric, positive definite matrix which implies the norm 
$ <u,v>_W= u^T W v$. 
From this  elementary follows 
\begin{equation}
<u,D_x v>_W = -<D_x^T u, v>_W  - u_1v_1 + u_N v_N  \label{SBPbyParts}
\end{equation}
which is the discrete analogous to partial integration.\footnote{A short derivation can be found in appendix~\ref{appendix:norm}.}
Setting $u = 1$ we find 
\begin{equation}
<1,D_x v>_W  =  -v_1 +v_N,    \label{SBPint}
\end{equation}
which is the discrete analogous of the fundamental theorem of calculus.
Thus, using SBP matrices, the discrete fluxes are well defined, depending only on the first or last point, and are in agreement with the analytical counterpart. 
By this, it follows that enforcing these boundary values allows to enforce the fluxes in and out of pipes.

 The SBP  concept shall be illustrated by a simple example.
Considering the discrete continuity equation 
\begin{equation}
    \partial_t \varrho + D_x (\varrho u) =0  , \label{mass1D}
\end{equation}
with $ \varrho , (\varrho u)  \in  \mathbb{R}^{N,1} $. 
To derive the conservation of the mass, equation \eqref{mass1D} is integrated, for which a discrete equivalent needs to be defined.  
The simplest choice, without considering SBP properties, is summing up, indicated by the vector ${\mathbf 1} \in\mathbb{R}^{N,1} $ with $ {\mathbf 1}_i = 1 \;  \forall i $,
\begin{equation}
    \partial_t {\mathbf 1}^T \varrho \Delta x  + \underbrace{\Delta x  {\mathbf 1}^T D_x}_{\tilde b^T} \varrho u =0  
\end{equation}
Choosing  $D_x$ as symmetric second order derivative with first order at the boundaries and $\Delta x$ the uniform grid spacing, this leads to
\begin{align}
    \tilde b^T =  
    \Delta x {\mathbf 1}^T D_x &= {\mathbf 1}^T \frac{\Delta x}{\Delta x} 
     \left(
        \begin{array}{ccccc} 
            -1 & 1 & \\
            -\frac{1}{2} & 0 & \frac{1}{2} & \\
            &-\frac{1}{2} & 0 & \frac{1}{2} & \\
            && \ddots & \ddots & \ddots   \\
            &&& -1 & 1   
        \end{array}
    \right) \\
    &= \left( -\frac 3 2 , \frac{1}{2} , 0 ,\cdots , 0, -\frac{1}{2} , \frac 3 2 \right). \label{eq:b_tilde}
\end{align}
The flux over the boundaries for this definition of the total mass $\tilde M = {\mathbf 1}^T \varrho \Delta x  $ yields $\partial_t \tilde M + {\tilde b^T} (\varrho u)  = 0 $, with 
\begin{gather}
{\tilde b^T} (\varrho u)  =  - \left( \frac 3 2(\varrho u )_1  - \frac 1 2 (\varrho u )_2   \right)
   +    \left( \frac 3 2(\varrho u )_{N}  -  \frac 1 2(\varrho u )_{N-1}   \right).
\end{gather}

For constant $\varrho u$, the fluxes at each boundary agree with the analytic fluxes. 
However, describing a wall only by  $u_1 = 0$ still  leaves a flux  of $\frac 1 2 (\varrho u )_2$, 
so that controlling the outermost points does not allow to control the mass flux into the pipe.    

To achieve zero flux for $u_1 = 0$ the derivative $D_x$ is decomposed  in terms of SBP to 
$ W^{-1} S= D_x$, with  
\begin{gather}
S = W D_x  =  \frac{1}{\Delta x} \left(
    \begin{array}{cccccc} 
    -\frac 1 2 & \frac 1 2 & \\
    -\frac 1 2  & 0 & \frac 1 2 & \\
    &-\frac 1 2 & 0 & \frac 1 2 & \\
    && \ddots & \ddots & \ddots &  \\
    &&& -\frac 1 2 & \frac 1 2   
    \end{array}\right) . \label{sbpM}
\end{gather}
The first and last lines are now only half a derivative, so that  
the weight matrix $W$ is 
\begin{gather}
  W_{ij}=
  \begin{cases} 
    &\frac 1 2 \text{ for } i=j=1, N \\ &1 \text{ for } i=j \neq 1, N \\ &0 \text{ for } i \neq j.
  \end{cases}
\end{gather}
This summation by part decompositions allows to rewrite the mass equation as 
\begin{gather}
    \partial_t W \varrho + S (\rho u ) = 0 . 
\end{gather}
Defining the mass as $ M = {\mathbf 1}^T W \varrho \Delta x  $  we obtain
 $\partial_t  M +  b^T (\varrho u)  = 0 $, with 
\begin{eqnarray}
b^T  (\rho u )= \Delta x \mathbf{1}^T S  (\rho u ) = (-1, 0, \dots,0, 1)  (\rho u ) =  -(\rho u )_1 + (\rho u )_N \, . 
\end{eqnarray}
Thus, the SBP matrices allow defining fluxes in and out of the pipes depending only on the values of the first and the last point. 
Discrete conservation laws can be derived where the fluxes determine the change of conserved quantities. 

To allow a decomposition like \eqref{SBPdecomp}, a special choice of the derivative operator is needed. Operators of different discretization error  are derived by Strand \cite{Strand1994} for diagonal matrices $W$, which are  preferred for our purpose,  since they allow an explicit prediction-correction scheme, see below. 
The derivation allows for arbitrary order, while we use second order in the numerical examples.
The discretization  in space is simply performed by replacing the derivatives $ \partial_x$ with the discrete  SBP derivative operators  $D_x = W^{-1}S$ in (\ref{eq_euler_conti}-\ref{eq_euler_momentum}).

The derived schemes are similar to a finite volume scheme. 
A subtle but essential difference is the double role of the outermost values which define the local discretization value and the flux, while for finite volumes the cell value and the flux at the boundary can be chosen independently. 

This allows a simple combination of systems, which can be illustrated by discretizing  the transport equation  on two one-dimensional domains or pipes ($\alpha=1,2$) governed by the transport equation 
$ \partial_t \varphi^\alpha + \lambda \partial_x \varphi^\alpha = 0   $, 
with transport velocity $\lambda$.  
It is discretized as    
\begin{equation}
	\partial_t W \varphi^\alpha + \lambda S \varphi^\alpha  = f^\alpha  . 
\end{equation}
A control flux $f^\alpha$ is introduced, which is nonzero only at the first and last point. 
It  is used to enforce the boundary conditions by choosing an appropriate value.  
We now  construct the connection of the  two domains using SBP matrices as before, see e.g.  \cite{NORDSTROM2009}.
This is done by assuming that the last point of the first domain and the first point of the second domain coincide, 
and further that the values of the two domains are identical in this point $\varphi_N^1\equiv \varphi_1^2 $.
We assume that $f_N^1 = - f_1^2 \equiv f $ so that the source in one domain cancels with the other to keep global conservation. 
This is later discussed in more detail. 
 The value of $f$ can be chosen to keep the value of $\partial_t \varphi^1_N $ and $\partial_t \varphi^2_1 $ identical. 
 This allows to combine the two equations as 
\begin{align}
\begin{split}
  \includegraphics[trim=0 40 0 0, clip,width=0.8\textwidth]{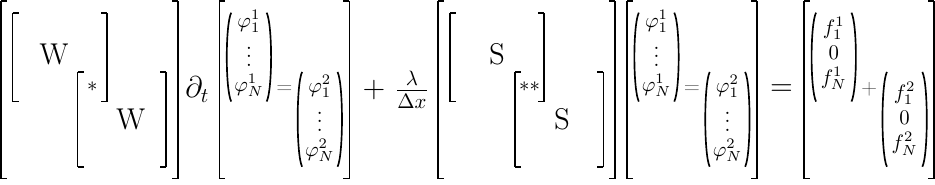}
\end{split}
\end{align}

\begin{equation*}
    *: W_{NN} + W_{11} = \frac{1}{2} + \frac{1}{2} = 1 \hspace{1ex}, \hspace{1ex} **:  S_{NN} + S_{11} = -\frac{1}{2} + \frac{1}{2} = 0 \hspace{1ex}, \hspace{1ex} f^1_N + f^2_1 =0
\end{equation*}
By this, effective matrices  for the total system can be formulated. 
For the second order central derivatives described above, the effective matrices are identical with a direct discretization as one domain. 
For other derivatives, the stencils might be modified close to the connection point.  
The essential fact is that with SBP schemes the domain can be controlled by controlling the outermost points, which can be done by setting the control fluxes. Setting these can  be done independently from solving the equations. 
This explains why no fluid problem needs to be solved to fulfill   the coupling condition. 
 For explicit derivatives, this can be done in a second step after evaluating the spatial part. 
 This is the core idea  for the construction of the gas pipe network scheme.
 Note that no smoothness was assumed for the coupling so that no  smoothness constraints in addition to the FD requirements are required. 

For the derivation of the gas pipes the discrete version of the delta function is needed. The definition of the integration 
including the diagonal weight matrix $W$ suggests to define 
\begin{equation}
\delta(x -x_j)_i \equiv  W^{-1}_{i,i}  \frac 1 {\Delta x}\delta_{i,j}    \label{discDelta}
\end{equation} 
which will be used below.

As a final remark, the property \eqref{SBPbyParts}   can be used to conserve quadratic invariants as the kinetic energy \cite{ReissSesterhenn2014}. 
This is not followed here since the total energy is not conserved due to the isothermal approximation. The integration by parts is, however, also the crucial step to derive the adjoint, so that discrete and analytical derivation become consistent.

\paragraph{Filtering}

The here used finite difference scheme as a direct discretization of the divergence form of the isothermal Euler equations is the simplest choice. 
It is augmented by a filter applied  between  time steps  to remove high-frequency oscillations, as commonly done in FD. 
The used central derivatives produce only small dissipation, also small anti-dissipation is possible. This is an effect of the non-linearity of the equations.  
The filter stabilizes the method and removes oscillations at locations with  a rapidly  varying behavior. 
In FV this is usually provided by the definition of the flux function, which inherently  creates the needed numerical dissipation. 
The filter approach of FD has the disadvantage that a sufficient  amount of dissipation needs to be adjusted by choosing the filter type or the filter  frequency, the advantage is that the negative definiteness of the filter  avoids  non-physical (negative)  friction and  entropy violating solutions are practically excluded.

For the network simulation, it is important that the filtering scheme does not change the values at the outermost points of all pipes. This is necessary to not alter the fluxes in the nodes.

We use two different filters, depending on the necessary filtering strength. For the smoother examples a Pade filter \cite{gaitonde2000pade} is applied. For the steeper gradient and the diamond network we use a conservative locally varying filter \cite{reiss2021pressuretight}:

\begin{equation}
    \overline{qA} = qA + \alpha_F D_L \cdot M(\sigma A) \cdot D_R \cdot q 
\end{equation}
with
\begin{equation}
    D_R = -D_L^T = \frac{1}{2}\begin{bmatrix} 
    -1 & 1 &  &  \\
     & \ddots& \ddots&\\
    &&-1&1
    \end{bmatrix} \hspace{1mm} , \hspace{2mm} M = \frac{1}{2}\begin{bmatrix} 
    1 & 1 &  &  \\
     & \ddots& \ddots&\\
    &&1&1
    \end{bmatrix} \hspace{2mm} \text{and} \hspace{2mm} \alpha_F = 0.5 \hspace{1mm} .
\end{equation}
The factor $\alpha_F$, $0<\alpha_F\leq 1$ is chosen big enough to avoid overshoots but small enough to be not too dissipative. With $\sigma$ the local filtering strength can be defined. To not influence the fluxes as necessary the first and the last two entries of that vector have to be set to zero. 
Setting only the outermost points of $\sigma $ zero would result in modified outermost points, but with a zero flux by  the filtering operation.
It can be viewed as an adiabatic boundary condition for the diffusion implied by the filter and might be useful in future implementations.

\subsection{Numerical implementation - pipes } 

The gas network is discretized using a finite-differences time-domain (FDTD) approach.
An explicit fourth-order Runge-Kutta scheme is employed for time-wise integration. 
However, also other time integration schemes are suitable.
The finite difference approach  used for the  spatial discretization is  based on the divergence form of the Euler equations (\ref{eq_euler_conti},\ref{eq_euler_momentum}). 
The derivatives $ \partial_x$ in (\ref{eq_euler_conti}-\ref{eq_euler_momentum}) are simply replaced with the discrete  SBP derivative operators  $D_x = W^{-1}S$.    
The utilized FD realization is not crucial for the following discussion and the following derivation  should be possible with many other FD forms.    
Essential is to  employ  a summation by parts scheme (SBP) for  the  discretization of the derivative operators, and further that this yields correct fluxes, as discussed below. 

The implementation of the coupling conditions (C1) and (C2) builds on additional fluxes non-zero only at the boundary points of the pipes, added to the right-hand-side of the discrete version of the governing Euler equations \eqref{eq_euler_conti}-\eqref{eq_euler_momentum} and $m=\varrho u$. 
\begin{align}
    \partial_t  A \varrho  &= - D_x A m + f_{\varrho} \nonumber \\ 
    \partial_t  A m        &= 
    -   D_x \left( A\frac{m^2}{\varrho}               \right)   
    - A D_x \left(                      c^2 \varrho \right) + f_m. \label{eq_euler_f} 
\end{align}
They might seem ad hoc, but they simply encode the information traveling in and out of the pipe. 
By this, these fluxes allow controlling the dynamics at the boundary, i.e. setting the boundary conditions.  We refer to them as control fluxes in the following.  
Since the numerical discretization of the pipes by using SBP operators provides well-defined fluxes,  analytical and numerical considerations are fully consistent, simplifying the argumentation. 
By this, the adjoint equations derived from the discrete scheme and the discretized adjoint equations obtained by discretizing the  analytical derived adjoint equations are consistent.
We assume in the following that the coupling conditions (C1) and (C2) are satisfied at the initial time.

As discussed above, for isothermal conditions the pressure and the density are directly linked by the speed of sound. 
Enforcing the same value for the density of all pipe endings in one junction forms the equal pressure condition (C1). 
To realize this constraint, the control  fluxes\footnote {The number of discretization points can vary between the pipes, however, we do not explicitly mark it ($N^\alpha\equiv N $).  } 
$f^\alpha_{\varrho,j}$, $j=1,N$ are non-zero  at only the first and last point of each pipe $\alpha$  of the continuity equation: 
\begin{equation}
\WW  \partial_t \left( \varrho A \right)_j^\alpha   =\WW  \underbrace{(- D_x \left( \varrho u A \right)^\alpha)_j}_{
  \widetilde{\partial_t \left( \varrho A \right)_j^{\alpha} } }
  + f_{\varrho,1}^{\alpha} \dfrac{\delta_{j, 1}}{\Delta x} 
  + f_{\varrho,N}^{\alpha} \dfrac{\delta_{j, N}}{\Delta x}. \label{eq_conti_with_add_fluxes}
\end{equation}
Therein, $\widetilde{\partial_t \left( \varrho A \right)_j^{\alpha}}$ denotes the unmodified right-hand-side which can be evaluated for each pipe separately.
Spatial integration over all pipes yields the total balance of mass which is changed by the physical fluxes and the control fluxes
\begin{equation}
 \partial_t \left( {\mathbf 1}^T   W   \left( \varrho A \right)^\alpha  \Delta x^\alpha \right) +
  \left[ \left( \varrho u A \right)^{\alpha} \right]_{x_{ 1}}^{x_{ N}}  = f_{\varrho,1}^{\alpha}  + f_{\varrho,N}^{\alpha} . \label{eq_c1_spatial_integration}
\end{equation}
Here, 
$ 
 {\mathbf 1}^T   W   \left( \varrho A \right)^\alpha  \Delta x^\alpha  
\equiv  
\int \left( \varrho A \right)^{\alpha} \mathrm d x 
 = M^\alpha 
$ is the mass in pipe $\alpha$. 

This form is essential for our scheme. 
The deployment of the SBP matrices created mass fluxes 
\begin{equation}
    \left[  \varrho u A \right]_{x_{ 1}}^{x_{ N}} = - ( \varrho u A)_{1}  + (\varrho u A)_{N}
\end{equation}
 depending,  in the discrete,  on the  values of the first and  last point only. 
 This allows controlling the fluxes  by controlling these points. 
Whatever FD scheme has this property could be used in the following, giving great freedom to adopt the scheme to special necessities.   
 

While the simple FD scheme  was found to work well,  
a more sophisticated   scheme could be used  if needed, as long as well defined fluxes are created by the SBP matrices. This includes many FD schemes, for example the skew symmetric scheme \cite{ReissSesterhenn2014} which has well-defined fluxes and avoids artificial dissipation by construction.

\subsection{Numerical implementation - junctions} 
The junctions are defined by the two coupling conditions C1 \eqref{equalPressure} and C2 \eqref{eq:conservation_mass}. These conditions are now enforced with the help of the control fluxes.  The well-defined fluxes guarantee that setting the outermost point controls the mass and momentum flux.      

 \paragraph{Coupling condition (C1) - equal pressure}

 The  mass  conservation in the full network gives an important restriction on the choice of the control fluxes.
The summation over all pipes of  the network gives the total mass conservation, as we assume no global in and out fluxes for the sake of simplicity. 
These could be included complicating the notation without changing the final result.
Thus, the change of total mass 
$M = \sum \limits_{{\alpha}} M^\alpha $, for a network of $N_p$ pipes ,  yields with \eqref{eq_c1_spatial_integration}
\begin{equation}
  \underbrace{\partial_t M}_{=0} +
  \sum \limits_{{\alpha=1}}^{N_p}  \sum \limits_{j=1,N}
  \left( \varrho u A \right)^{\alpha}_j n_j 
  - 
  \sum \limits_{{\alpha=1}}^{N_p}  \sum \limits_{j=1,N} f_{\varrho,j}^{\alpha} = 0
\end{equation}
with $[ (\varrho u A)_j^\alpha ]_{x_{j=1}}^{x_{j=N}} = \sum_{j=1,N} \left( \varrho u A \right)^{\alpha}_j  n_j  $, $n_1 = -1 $ and $n_N = 1 $.  
Reordering the terms to form sums over each junction $k=1,\dots, K$ yields\footnote{Assuming all pipes belong to at least one node.}
\begin{equation}
    \sum \limits_{k=1}^K \sum \limits_{(\alpha,j) \in N_k} \left( \varrho u A \right)^{\alpha}_j n_j 
    - \sum \limits_{k=1}^K \sum_{(\alpha,j) \in N_k} f_{\varrho,j}^{\alpha}  = 0.
\end{equation}
The first term is zero due to Kirchhoff's law (\ref{Kirchhoff}), so that  
\begin{equation}
  \sum \limits_{k=1}^K \sum \limits_{(\alpha,j) \in N_k} f_{\varrho,j}^{\alpha}  = 0
\end{equation}
follows.
 With the additional assumption that the mass conservation is fulfilled for each node separately, excluding that a mass defect in one node is compensated by another node,  
\begin{equation}
  \sum \limits_{(\alpha,j) \in N_k} f_{\varrho,j}^{\alpha} = 0 \label{eq_node_wise_kirchhoff}
\end{equation}
is found.
Thus, the sum over all control fluxes at a junction has to vanish. 
This zero-sum of the control fluxes was already used in the introductory  example of the transport equation.
\smallskip 

Now, we derive the mass control fluxes per node. 
For this, we consider the start or end of a pipe ($j=1$ or $j = N$) belonging to a node of interest by which we keep only one of the control fluxes:
\begin{equation}
W_{jj}     \partial_t \left( \varrho A \right)^{\alpha}_j = 
W_{jj}     \widetilde{\partial_t \left( \varrho A \right)^{\alpha}_j} 
     + 
     \dfrac{f^{\alpha}_{\varrho,j}}{\Delta x^\alpha} \label{eq_conti_with_flux}
\end{equation}
Corresponding to \eqref{equalPressure} all pipes in a junction $N_k$ should have the same time derivative 
$ \partial_t \left( \bar \varrho \right)_k $  
\begin{equation}
  \partial_t \left( \varrho \right)_j^{\alpha}\left.\right|_{(\alpha,j) \in N_k} \equiv \partial_t \left( \bar \varrho \right)_k. 
\end{equation}
Using \eqref{eq_conti_with_flux} and assuming a time independent  cross sections $A$  yields
\begin{align}
   W_{jj} A^\alpha \partial_t \left( \bar \varrho \right)_k 
  = 
  W_{jj}  \widetilde{\partial_t \left( \varrho  A \right)_j^{\alpha}}\left.\right|_{(\alpha,j) \in N_k} + \dfrac{f^{\alpha}_{\varrho,j}}{(\Delta x)^{\alpha}}\left.\right|_{(\alpha,j) \in N_k}.   
\end{align}
Summing over all pipe endings, which  belong to the junction $N_k$ and multiplication with the grid spacing results in 
  \begin{align} 
  \sum \limits_{(\alpha,j) \in N_k}  W_{jj}  A^{\alpha}\Delta x^{\alpha}\partial_t \left( \bar \varrho \right)_k 
  = 
  \sum \limits_{(\alpha,j) \in N_k}  W_{jj}  \Delta x^{\alpha}\widetilde{\partial_t \left( \varrho  A \right)_j^{\alpha}}+ \underbrace{\sum \limits_{(\alpha,j) \in N_k}f^{\alpha}_{\varrho,j} }_{=0} 
\end{align}
wherein \eqref{eq_node_wise_kirchhoff} was used. 
Thus, the average density derivative results, with $W_{1,1} = W_{N,N}$, to
\begin{equation}
\partial_t  \left( \bar \varrho \right)_k 
= 
\dfrac{\sum \limits_{(\alpha,j) \in N_k}\Delta x^{\alpha}\widetilde{\partial_t \left( \varrho A \right)^{\alpha}_j}}{\sum \limits_{(\alpha,j) \in N_k}\Delta x^{\alpha} A^{\alpha} } .
\end{equation}
This can be used to calculate the fluxes 
\begin{equation}
     f_{\varrho,j}^{\alpha}   
 = 
\WW_{jj} {\Delta x^\alpha}  \left( A^\alpha \dfrac{\sum \limits_{(\alpha,j) \in N_k}\Delta x^{\alpha}\widetilde{\partial_t \left( \varrho A \right)^{\alpha}_j}}{\sum \limits_{(\alpha,j) \in N_k} A^{\alpha}\Delta x^{\alpha} }
   -   
  \widetilde{\partial_t \left( \varrho A \right)_j^{\alpha} } ) \right)
\label{fluxRho}
\end{equation}
or, yielding the same result, can be used to overwrite the right hand side $\widetilde{\partial_t \left( \varrho \right)_j^\alpha}$ directly resulting in the constraint (C1).

\paragraph{Coupling condition (C2) -  conservation of mass}

To realize the coupling condition (C2) the control fluxes $f_m$ of  the momentum equation \eqref{eq_euler_f}  are used as for the mass equation,
\begin{equation}
  W_{jj} \partial_t \left( \varrho u A \right)^{\alpha}_j
  =
  W_{jj} \widetilde{\partial_t \left( \varrho u A  \right)^{\alpha}_j} + \dfrac{f^{\alpha}_{m,j}}{\Delta x^{\alpha}} \label{eq_momentum_with_flux}
\end{equation}
Again, $
 \widetilde{\partial_t \left( \varrho u A  \right)^{\alpha}_j } = 
\left( 
    -   D_x \left( A\frac{m^2}{\varrho}               \right)   
    - A D_x \left(                      c^2 \varrho \right) 
\right)_j     $ 
corresponds to the (spatial) right-hand-side without coupling conditions.
Summing over all pipe endings belonging to a junction $N_k$, as above in \eqref{eq_c1_spatial_integration}, results in
\begin{equation}
    \sum \limits_{(\alpha,j) \in N_k} 
     \left ( 
              W_{jj}  {\partial_t \left( \varrho u A\right)^{\alpha}_j n_j} 
      \right)
    = 
    \sum \limits_{(\alpha,j) \in N_k} \left ( W_{jj} \widetilde {\partial_t \left( \varrho u A\right)^{\alpha}_j n_j} + \dfrac{f^{\alpha}_{m,j}}{\Delta x^{\alpha}} \right)
    \equiv 0.
\end{equation}
The whole expression is zero, since the left-hand-side corresponds to the time derivative of Kirchhoff's law \eqref{Kirchhoff}, as all
${W_{jj},~j=1,N}$ are equal at the end points of the pipes. 
Separating the flux terms leads to
\begin{equation}
  \sum \limits_{(\alpha,j) \in N_k} W_{jj}\widetilde {\partial_t \left( \varrho u A\right)^{\alpha}_j n_j} 
  + 
  \sum \limits_{(\alpha,j) \in N_k} \dfrac{f^{\alpha}_{m,j}}{\Delta x^{\alpha}}  = 0.
\end{equation}
By assuming equal $f^\alpha_{m,j} = \bar f_m$ for in- and out-fluxes
\begin{equation}
  \sum \limits_{(\alpha,j) \in N_k} W_{jj}\widetilde {\partial_t \left( \varrho u A\right)^{\alpha}_j n_j} + \sum \limits_{(\alpha,j) \in N_k} \dfrac{\bar f_m}{\Delta x^{\alpha}}   = 0
\end{equation}
is found.
However, this relation is only one possible choice. 
Physically it is clear, that if we find a mass defect at a node changing any of the mass fluxes of the connected pipes (or all together) can correct this.\footnote{Formally, the mass conservation (C2) is one condition, while the equal pressure condition is in fact 
number of pipes in one node minus one condition. }
Real nodes might behave differently, in particular depending on the detailed geometry of the junction, an aspect which is often not modeled in networks. 
However, when testing different distributions, we found the simulation results surprisingly insensitive to this choice.
For our choice the  control flux results to
\begin{equation}
  \bar f_m = - \frac{\sum \limits_{(\alpha,j) \in N_k}  W_{jj} \widetilde {\partial_t \left( \varrho u A\right)^{\alpha}_j n_j}}{\sum \limits_{(\alpha,j) \in N_k} \dfrac{1}{\Delta x^{\alpha}_j} } .
\label{fluxM}
\end{equation}
In conclusion, the old right-hand-sides are modified by $\dfrac{1}{\Delta x^{\alpha}_j}\bar f_m $ at each pipe in node $k$ to ensure the coupling condition (C2).

\paragraph{Coupling conditions - linear coupling matrix}

\begin{figure}
\begin{center}
\includegraphics[width=0.5\textwidth]{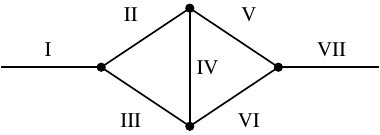}
\caption{Example: diamond}
\label{fig:diamond}
\end{center}
\end{figure}

\begin{figure}[bh]
    \centering
    \includegraphics[trim = 13cm 1.5cm 12cm 1.5cm, clip, width=0.5\textwidth]{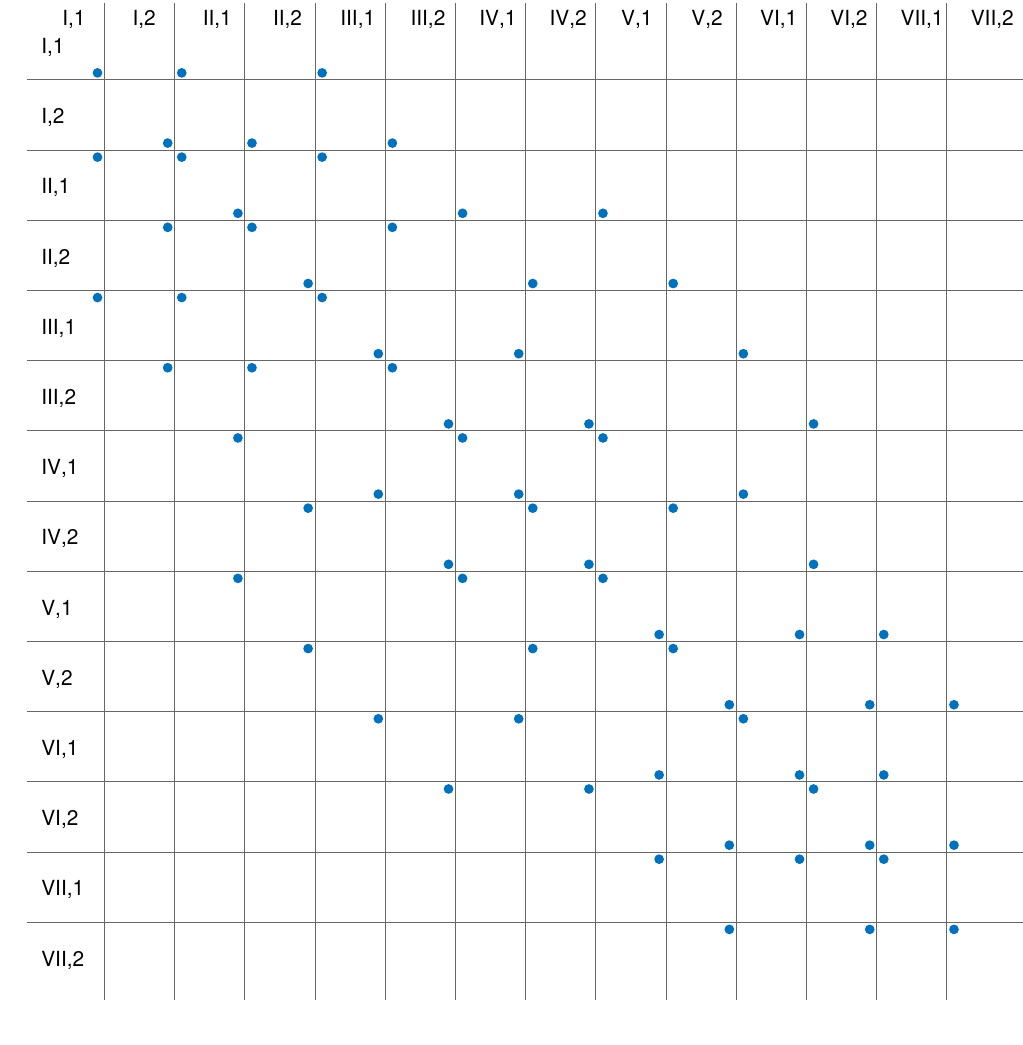}
    \caption{Sparse coupling matrix C for the diamond network. The dots show the nonzero elements. The roman numbers refer to the pipe numbers, 1 and 2 to the variables $\varrho$ and $\varrho u$.}
    \label{fig:coupling}
\end{figure}

The above findings enable a simple predictor-corrector scheme for the coupling of pipes.
Since the fluxes are linear in the undisturbed right-hand-side, the two conditions are realized in additional control fluxes $f$ which can be rewritten in a  matrix $C$ that contains the coupling information. 
\begin{equation}
   \dot{Q}  = \widetilde{\dot{Q}} + F   = \widetilde{\dot{Q}} + C \cdot \widetilde{\dot{Q}} . \label{eq:qdot_network}
\end{equation}
Therein, $\widetilde{\dot Q}$ is the right-hand-side without the control  fluxes with 
$ \widetilde {\partial_t \left( \varrho  A\right)^{\alpha} }$  
and
$ \widetilde {\partial_t \left( \varrho u A\right)^{\alpha} }$ for the whole network in one vector.
$\dot Q$ holds the whole information of the system with the control fluxes, i.e. coupling conditions.  
Since the control fluxes are linear in $ \widetilde{\dot{Q}}$ , see \eqref{fluxRho} and \eqref{fluxM},  their effect can be rewritten in a (sparse) matrix $C$. It is shown for the diamond network in figure \ref{fig:coupling} which will be discussed in section \ref{sec:E2}. 
The numbering of the pipes is according to figure \ref{fig:diamond}.
Taking one of the labels $I,2$ for example means that this part of the matrix refers to the first pipe but the second equation of \eqref{eq_euler_f} ($\partial_t {\varrho u}$). 
The  entries of the matrix are either on the first or the last position of the block belonging to one pipe and variable. 
For example the first entry of $V,1$ is connected to the last entry of $II,1$. 
The value of the entry depends on the number of pipes in the node and the area of the pipe. 

To simulate the network at first the right-hand-sides of the Euler equations \eqref{eq_euler_conti}-\eqref{eq_euler_momentum} are calculated for the individual pipes. 
The values are stored in one vector $\widetilde{\dot{Q}}$. 
By the coupling matrix $C$ the vector  $\dot{Q}$ is calculated, which is used for the time-stepping scheme.

With that, the code is well suited for parallelization because all pipes can be solved independently, and only  after this  a single, sparse exchange is needed. 
Furthermore, this form can be used to derive the adjoint equations compactly.

\section{Adjoint approach}

Adjoint equations allow calculating how a target, described by a scalar function of the system state, changes with a high dimensional modification of the system equation.  
This is helpful in design and control since it shows how to approach a given control or design goal.   
In general, adjoint equations can be derived in different ways, the continuous or the discrete approach. Also, automatic differentiation techniques are used to derive adjoint codes. Despite different discretizations, all approaches are consistent and applicable, see \cite{GilesPierce2000}.
  

\paragraph{General idea}
Here, the adjoint equations are introduced in discrete version \cite{GilesPierce2000} using a matrix-vector notation.

Adjoint equations arise by a scalar-value objective given by the product between a geometric weight $g$ and the system state $q$:
\begin{equation}
  J = g^\text{T}q, \qquad g,q\in\mathbb{R}^{n}
\end{equation}
The system state $q$ is the solution of the governing system
\begin{equation}
  L q = r, \qquad L\in\mathbb{R}^{n \times n}, ~~r\in\mathbb{R}^{n}
\label{eq_linear_system}
\end{equation}
with $L$ being the governing operator e.g.~the (linearized) isothermal Euler equations. 
In this case, the discrete solution vector in space and time is written in one vector. 
A source term $r$  is added  on the right-hand-side as a modification of the governing equation.
The corresponding adjoint is defined by
\begin{equation}
  L^\text{T}q^* = g,
	\label{eq_adjoint_system_1}
\end{equation}
with the adjoint variable $q^*$.
Using the objective and the governing system equation with a Lagrangian multiplier leads to
\begin{equation}
  J = g^\text{T}q 
  =
  g^\text{T}q 
  - {q^*}^T\left(  L q - r
\right)
=  q^T\underbrace{( g - L^T q^*  )}_{\equiv 0 }    
+ {q^*}^\text{T}r 
= {q^*}^\text{T}r .
\label{eq_adjoint_simple}
\end{equation}
With that, the objective function $J$ can be computed for every possible $r$ once the adjoint equation is solved. 
So gradients for $J$ with respect to $r$ can be computed efficiently with $q^*$.


\paragraph{Objective}
In the following, the objective function $J$ under consideration is defined in space and
time with $\mathrm d \Omega = \mathrm d x \mathrm d t$ over all pipes:
\begin{equation}
  J = \frac{1}{2} \iint \left( q - q_{\mathrm{target}} \right)^2 \sigma \mathrm
d \Omega. \label{eq_objective}
\end{equation}
The additional weight $\sigma(x,t)$ defines where and when the objective is evaluated.
The goal  could be keeping a certain reference pressure, but more involved objective functions are possible. 
We use a spatial discrete vector notation  so that the former expression becomes 
\begin{equation}
  J = \frac{1}{2} \int \left( q - q_{\mathrm{target}} \right)^T 
  (W \sigma )
  \left( q - q_{\mathrm{target}} \right)
  \mathrm
d t. \label{eq_objective_W}
\end{equation}
where $(W \sigma) $ is a diagonal operator, with the according weights.  
The scalar product is written in matrix notation where the index runs over the spatial discretization points and the variables $\rho A$ and $Am= A\rho u$. 
The variable $q_{\mathrm{target}}$ denotes a desired system state, which is to be reached by optimal modification of sources $r$.

In practice, the objective function is supplemented by additional constraints, e.g.~control effort, and a regularization term \cite{Bewley2001}. 
For the latter, the control variable $r$ is added to the objective function using a weight parameter $l$
\begin{equation}
  J = \ldots + \dfrac{l^2}{2} \iint r^2  ~\mathrm d \Omega.
\end{equation}
For the examples presented below, $l$ is chosen to a suitable, small value.

An optimal system state is realized if $J$ reaches a minimum.
The minimum is to be achieved under the constraint
that the isothermal Euler equations \eqref{eq_euler_conti}-\eqref{eq_euler_momentum} are satisfied.


\paragraph{Adjoint}
For the concrete derivation of the adjoint equations, the isothermal Euler equations \eqref{eq_euler_conti}-\eqref{eq_euler_momentum} are amended by a force $r(x,t)$ on the right-hand-side and  abbreviated with $\mathcal{F}(Q)$. 

\begin{equation}
  \partial_t q +  D_x \mathcal{F} (q)   
  = r  \label{eq:euler_start_adjoint}
\end{equation}
with $q = [A \varrho, A m]$, 
$ \mathcal{F} = [A m, A m^2/\varrho + c^2 \varrho A]$ 
and 
$r = [r_\varrho, r_m]$.
The linerization with respect to $q + \delta q$ yields
\begin{equation}
W    \partial_t \dfrac{\partial q_i}{\partial q_j}\delta q_j + S \dfrac{\partial \mathcal{F}_i}{\partial q_j}\delta q_j = W \delta r,
\end{equation}
with $W$ and $S$ being the matrices introduced in section \ref{sec:numerical_implementation} results in
\begin{align}
 \pmb{\mathcal{W}}   \partial_t \underbrace{\begin{pmatrix} 1 & 0 \\ 0 & 1 \end{pmatrix}}_{L_1} \delta q 
 + 
\pmb{\mathcal{S}} \underbrace{\begin{pmatrix} 0 & 1 \\ -\frac{m^2}{\varrho^2} + c^2 & \dfrac{2m}{\varrho} \end{pmatrix} }_{L_2}\delta q &= \pmb{\mathcal{W}}\delta r 
\end{align}
or
\begin{align}
\pmb{\mathcal{W}}    \partial_t L_1 \delta q + \pmb{\mathcal{S}} L_2 \delta q &= \pmb{\mathcal{W}} \delta r \label{eq:euler_linearized}
\end{align}
with the block matrices 
\begin{align}
 \pmb{\mathcal{W}}   = \begin{pmatrix} W & 0 \\ 0 & W \end{pmatrix} , 
 \; \; 
\pmb{\mathcal{S}} = 
\begin{pmatrix} S & 0  \\ 0&  S \end{pmatrix} , 
\; \; 
\pmb{\mathcal{B}} = 
\begin{pmatrix} B & 0  \\ 0&  B \end{pmatrix} .
\label{SandW}
\end{align}
The last matrix is defined for later use. 

These equations are combined with the linearized objective \eqref{eq_objective}, given by
\begin{equation}
    \delta J = \dfrac{\partial J}{\partial q}\delta q = \int \underbrace{(q-q_{\mathrm{target}})^T\sigma}_{g^T} \pmb{\mathcal{W}}\delta q \mathrm{d}t\label{eq:lin_objective}
\end{equation}
in a Lagrangian manner, using a Lagrangian multiplier $q^* = [\varrho^*, m^*]$, which becomes the adjoint variable.
\begin{equation}
    \delta J = \int g^T \pmb{\mathcal{W}}\delta q \mathrm{d} t 
    - 
    \int {q^*}^T
                \left( \pmb{\mathcal{W}} \partial_t \delta q 
                     + \pmb{\mathcal{S}}  L_2 \delta q 
                     - \pmb{\mathcal{W}}\delta r \right) \,\mathrm{d}t \label{eq:lagrangian}
\end{equation}
In order to remove the dependency on $\delta q$ we aim to reorder the terms.  
Since the integrand is a scalar it can be formally transposed without changing it. 
Using standard rules for the transpose operation and by using \eqref{SBPtrans} and \eqref{SandW} to get 
 $\pmb{\mathcal{S}}^T = -\pmb{\mathcal{S}} + \pmb{\mathcal{B}} $ we find       
\begin{eqnarray}
    \delta J &=& 
    \int \delta q^T \left( \pmb{\mathcal{W}} g + \pmb{\mathcal{W}} \partial_t q^* + L_2^T \pmb{\mathcal{S}} q^* \right) \mathrm{d}t  + \int {q^*}^T\pmb{\mathcal{W}} \delta r\, \mathrm{d} t  \no\\
    &-& \left[ {q^*}^T \pmb{\mathcal{W}}\delta q \right]_{t_0}^{t_{\mathrm{end}}} 
    - 
    \int  {q^*}^T \pmb{\mathcal{B}} L_2 \delta q  \mathrm{d} t
\end{eqnarray}
by partial integration.%
The dependency of $\delta q$ is removed by demanding
\begin{equation}
    \pmb{\mathcal{W}}  g + \pmb{\mathcal{W}}  \partial_t q^* + L_2^T \pmb{\mathcal{S}}   q^* \stackrel{!}{=} 0
\end{equation}
resulting in the adjoint isothermal Euler equations\footnote{With $W^{-1}L_2^T W = L_2^T$ - this holds for our diagonal matrix $W$}
\begin{equation}
    \partial_t q^* = - L_2^T D_x q^* - g.
    \label{eq_euler_adjoint_matrix} 
\end{equation}

The integrals resulting from partial integration have to vanish as well. This gives rise to the adjoint boundary and initial conditions.
Adjoint boundary conditions are discussed elsewhere \cite{Lemke2015}. In the appendix, the adjoint boundary conditions for a non-reflecting boundary are provided.
The temporal (initial) condition of the adjoint system is given at the final time. In general, the adjoint system is well-posed only if the adjoint initial state
is defined at the end of the computational time and the system is integrated backwards \cite{GilesPierce2000}.

The resulting adjoint equations for the isothermal Euler equations are in analytical form 
\begin{align}
    \partial_t  \varrho^*   &= \left(u^2 - c^2\right) \partial_x m^* - g_\varrho  \nonumber \\ 
    \partial_t  m^*         &= -\partial_x \varrho^* -2u \partial_x m^* - g_m  \label{eq_euler_adjoint} 
\end{align}

$g_\varrho$ and $g_m$ are the linearized parts of the objective function $J$ \eqref{eq:lin_objective}. Depending on $J$ they can occur in both equations or one of them is zero.

\paragraph{Adjoint network}

For the adjoint network we start again with equation~\eqref{eq:qdot_network} and replace $\widetilde{\dot{Q}}$ with \eqref{eq:euler_start_adjoint} but with $Q$ for the whole network and $\pmb{\mathcal{D}_x}$ another block matrix:
\begin{equation}
    \dot{Q} = -\pmb{\mathcal{D}_x} \mathcal{F}(Q) + R + C \cdot \left(-\pmb{\mathcal{D}_x} \mathcal{F}(Q) + R\right) . 
\end{equation} 
As Q, the capital R collects the source terms  of the whole  network. Using SBP matrices this leads to
\begin{equation}
    \pmb{\mathcal{W}} \dot{Q} = (1+C)(-\pmb{\mathcal{S}} \mathcal{F}(Q) + \pmb{\mathcal{W}} R) 
\end{equation} 
$\pmb{\mathcal{W}}$ and $\pmb{\mathcal{S}}$ are still block matrices but with accordingly increased  number of blocks to fit the network. Linearization leads to
\begin{equation}
    \pmb{\mathcal{W}}    \partial_t \pmb{\mathcal{L}_1} \delta Q =
    (1+ C) \left( - \pmb{\mathcal{S}} \pmb{\mathcal{L}_2} \delta Q +  \pmb{\mathcal{W}} \delta R \right)
\end{equation}
with $\pmb{\mathcal{L}_1}$ and $\pmb{\mathcal{L}_2}$ being block matrices with $L_1$ and $L_2$ respectively as defined in equation \eqref{eq:euler_linearized}.

Using the linearization, the objective and a Lagrangian multiplier $Q^*$ yields:
\begin{equation}
    \delta J = \int G^T \pmb{\mathcal{W}} \delta Q \mathrm{d} t - \int {Q^*}^T\left( \pmb{\mathcal{W}} \partial_t \pmb{\mathcal{L}_1} \delta Q +
    (1+C) (\pmb{\mathcal{S}} \pmb{\mathcal{L}_2} \delta Q -  \pmb{\mathcal{W}} \delta R) \right) \mathrm{d}t 
\end{equation}
The capital letters $G$, $Q$ and $R$ refer again to the whole network. Partial integration with neglecting the boundaries leads to:
\begin{align}
    \delta J &= \int \delta Q^T \left(\pmb{\mathcal{W}^T} G + \pmb{\mathcal{L}_1^T} \partial_t \pmb{\mathcal{W}^T} Q^* +
    \pmb{\mathcal{L}_2^T} \pmb{\mathcal{S}^T} (1+C^T) Q^* \right) \mathrm{d} t
    \no\\&\qquad + \int {Q^*}^T \left( \pmb{\mathcal{W}} \delta R + C \pmb{\mathcal{W}} \delta R \right) \mathrm{d} t . 
    \end{align}
    By demanding the adjoint equation to be zero, the dependency on $\delta Q$ is removed.
    \begin{align}
    \partial_t \pmb{\mathcal{W}} Q^* &= - \pmb{\mathcal{W}} G -  \pmb{\mathcal{L}_2^T} \pmb{\mathcal{S}^T} \left( 1 + C^T  \right)Q^*
\end{align}
This is the same as applying the adjoint Euler equations~\eqref{eq_euler_adjoint_matrix} on $Q^* + C^T \cdot Q^*$ instead of $q^*$.  

\begin{figure}
  \centering
  \includegraphics[width=.75\textwidth]{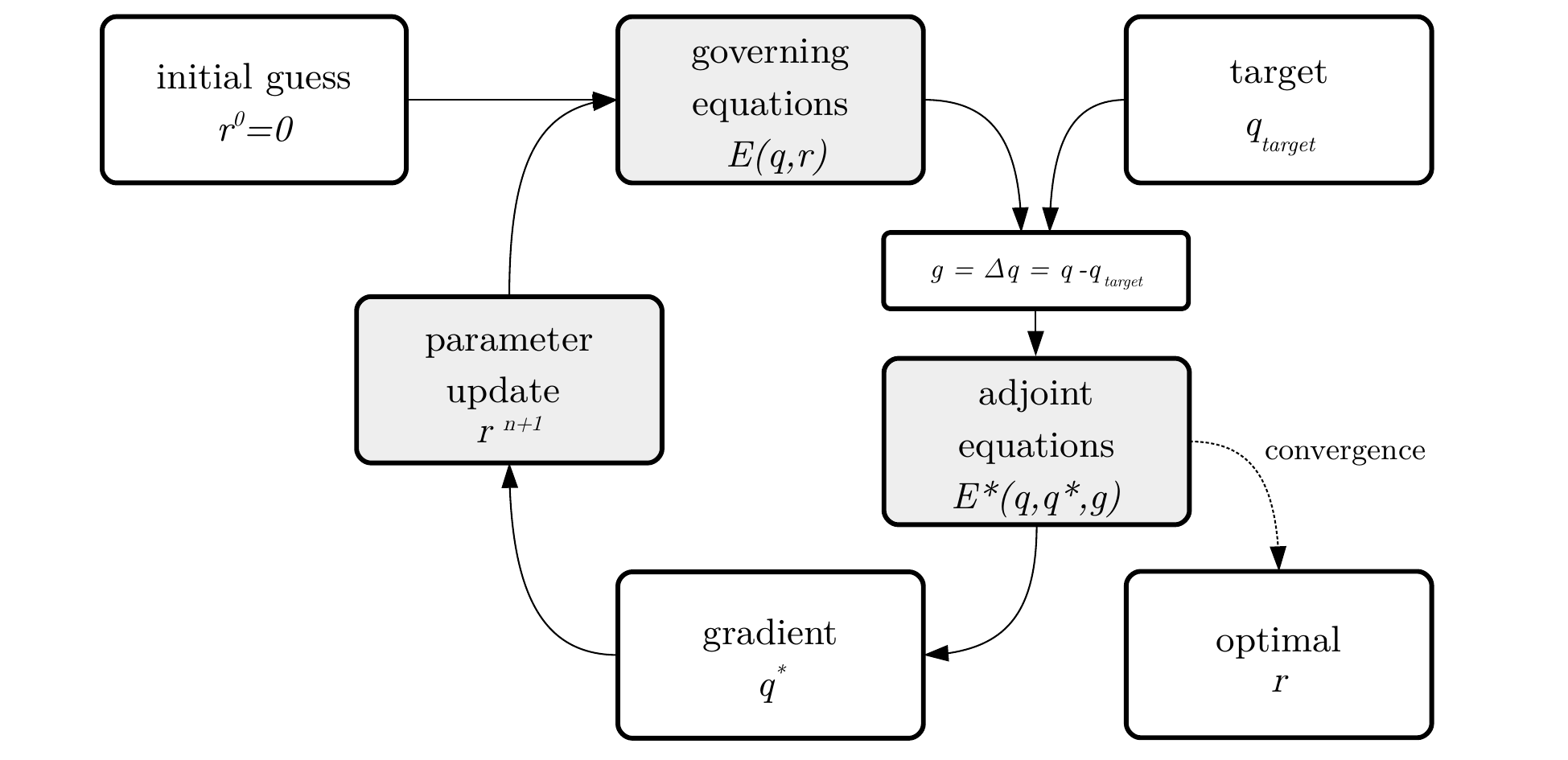}
  \caption{Iterative data assimilation procedure. 
           Expensive operations are marked by a grey box. 
           Details are given in the text. \label{fig_iterative_framework}}
\end{figure}

\paragraph{Iterative procedure}
The adjoint is a high dimensional gradient. 
In principle, any gradient based optimization method can be used to obtain the optimum of $J$. 
The most basic method is the steepest descent, where the direction of the gradient is directly used iteratively. 
First, the governing isothermal Euler equations \eqref{eq_euler_conti}-\eqref{eq_euler_momentum} are solved forward in time taking into account all constraints.
Subsequently, the adjoint equations \eqref{eq_euler_adjoint} are calculated backward in time incorporating the direct solution and the weight $g$ resulting from the considered objective \eqref{eq_objective}.
Based on the adjoint solution, the gradient $\nabla_r J$ is determined and used to update $r_{\varrho / m}^n$:
\begin{equation}
  r_{\varrho / m}^{n+1} = r_{\varrho / m}^n + \alpha_s \nabla_r J \theta_I,	\label{eq_steepest_descent}
\end{equation}
with $\alpha_s$ denoting an appropriate step size and $n$ the iteration number.
The gradient is calculated for the whole computational domain and the full
simulation time, but only evaluated at $\theta_I$ the with "I" labeled areas in figure~\ref{ex1_setup} and \ref{fig:E2_setup}.
The procedure is repeated until convergence is reached.
The iterative procedure is illustrated in Fig.~\ref{fig_iterative_framework}.

Please note, the method optimizes towards local extrema.
Determining a global optimum is not ensured.
The computational costs of the adjoint approach are independent of the number of parameters to be optimized but on the size of the computational domain and the number of time steps to be carried out.

For more information on the boundary conditions and the derivation of the adjoint equations for the (non isothermal) Euler and Navier-Stokes equations see Lemke\cite{Lemke2015}.

\section{Examples}

In the following, three examples are presented to show the applicability of the introduced finite difference technique.
The adjoint approach is used for simple but representative optimization tasks.
A Pade filter  or a conservative local varying filter is used after each  time step to guarantee a smooth solution, as discussed above. The CFL number is chosen to $0.77$, see table~\ref{tab:parameter}.

\subsection{(E1.1) Three pipes in a row}
\begin{figure}[ht]
\centerline{\includegraphics[width=0.7\textwidth]{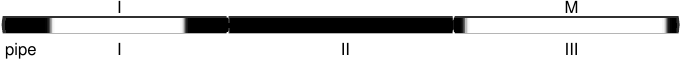}}
\caption{Setup 3 pipes: I: influence of the adjoint, M: measurement, domain of the objective function $J$}
\label{ex1_setup}
\end{figure}

The first example considers three pipes in a row, see figure~\ref{ex1_setup}. 
The goal is to create an acoustic density pulse in pipe \RM{3}  by a force $r$ in pipe \RM{1}.
The objective function is defined by
\begin{equation}
    J = \dfrac{1}{2} \iint \left( \varrho - \varrho_{target}\right)^2 \sigma ~ \mathrm d \Omega. \label{eq_examples_objective}
\end{equation}
Therein, $\varrho_{target}(x) = \varrho_{ref} \cdot (1+\beta \cdot exp( -(x-x_0)^2/\gamma^2))$, denotes the target density profile with $x_0  = 5\pi$, $\beta= 0.03$ and $\gamma = 0.9$, $\sigma(x,t)$ a weight defining where and when the objective is to be evaluated and $\Omega$ the space-time-measure for the whole computational domain and time.

The weight $\sigma$ is chosen to be non-zero only for the penultimate time step.
Spatially it is defined as tophat with a smooth fade-in-fade-out in pipe \RM{3}, see (M: measurement) in figure~\ref{ex1_setup}.

The adjoint-based gradient \eqref{eq_steepest_descent} 
\begin{equation}
    r^{n+1} = r^{n} + \alpha_s \nabla_r J \theta
\end{equation}
is evaluated in terms of the steepest descent approach with a suitable step size $\alpha_s$.
The weight $\theta(x)$ controls the location of the forcing which is restricted to pipe \RM{1}.
The corresponding region (I: influence) with a smooth fade-in-fade-out is shown in figure~\ref{ex1_setup}.
The main parameters of the simulation are shown in table~\ref{tab:parameter}.

\begin{center}
\begin{table}[t]%
\centering
\caption{Parameters for the examples (E1) and (E2). 
$N$ and $L$ are the number of grid points and length for each pipe.\label{tab:parameter}}%
\begin{tabular}{lccccc}
\toprule
\textbf{Case} & \textbf{time steps} & \textbf{cfl}  & \textbf{N}  & \textbf{L}  & \textbf{$\alpha_s$ } \\
\midrule
E1.1: 3 pipes & 1000 &$0.77$ & 200 & $2\pi$ & $3\cdot 10^6$ \\
E1.2: 3 pipes, steep gradient & 1000 &$0.77$ & 200 & $2\pi$ & $3\cdot 10^2$ \\
E2: diamond & 5000 &$0.77$ & 200 & $2\pi$ & $2.5\cdot 10^4$ \\
E2: diamond & 10000 &$0.77$ & 400 & $2\pi$ & $2.5\cdot 10^4$ \\
E2: diamond & 20000 &$0.77$ & 800 & $2\pi$ & $2.5\cdot 10^4$ \\
\bottomrule
\end{tabular}
\end{table}
\end{center}


\begin{figure}
    \includegraphics[width = 0.49\textwidth]{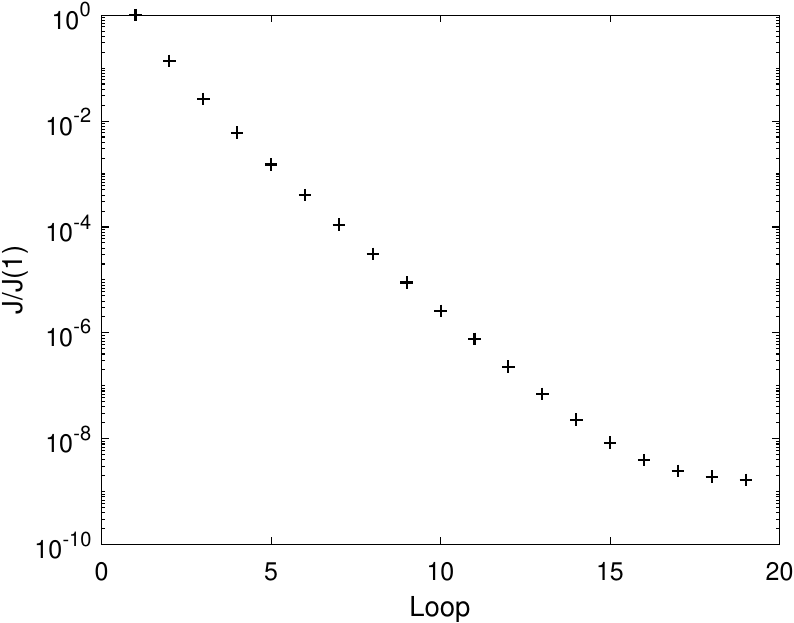}
    \includegraphics[width = 0.49\textwidth]{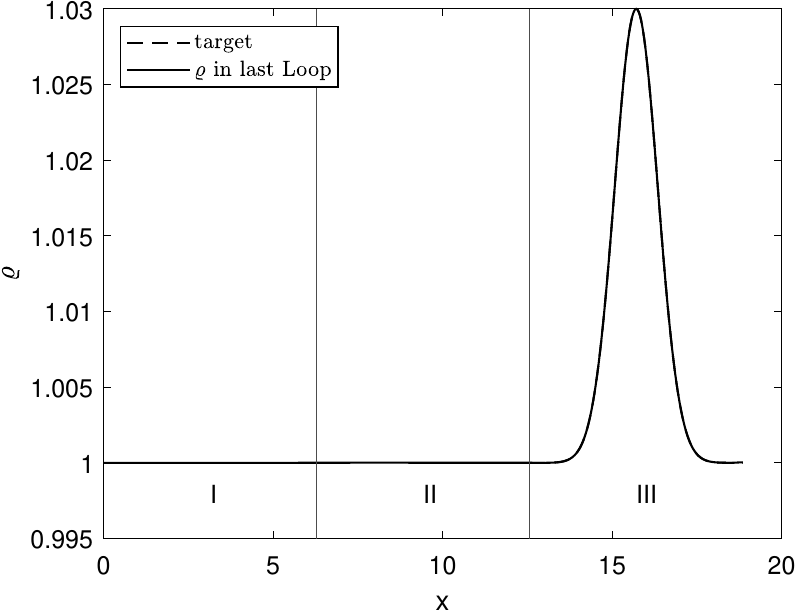}
    \caption{(Left) Progress of the objective function of (E1) normalized with respect to initial loop $J/J_1$. (Right) Density profile $\varrho$ at the penultimate time step after 19 iterative loops.
    \label{fig:E1_J}
    }
\end{figure}

After 19 iterative loops of the adjoint-based framework, convergence is achieved in terms of the objective function $J$ using the criterion 
\begin{equation}
    \frac{J_{n-1}-J_n}{J_1} < 2\cdot 10^{-10} . 
\end{equation}
The objective is reduced by nearly 9 orders of magnitude, see figure~\ref{fig:E1_J} (left). 
The desired pressure distribution occurs in pipe~\RM{3} in the penultimate time step, see figure~\ref{fig:E1_J} (right).
The remaining derivation is in the order of $10^{-6}$ kg/m$^3$.
Accordingly, the optimization target was achieved.

\begin{figure}
    \includegraphics[width = 0.49\textwidth]{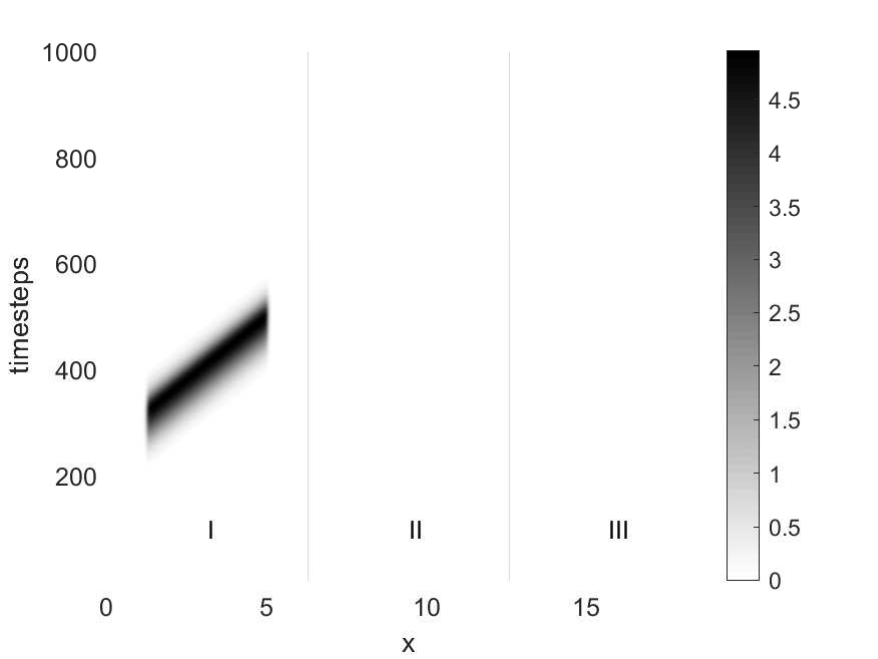}
    \includegraphics[width = 0.49\textwidth]{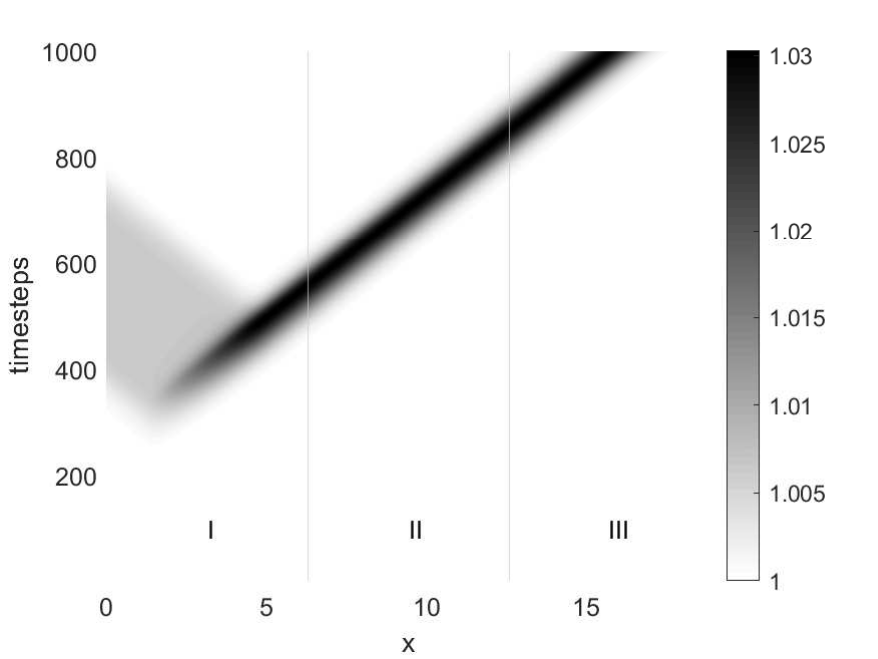}
    \caption{(Left) (E1) $x-t-$diagram for the excitation $r$. (Right) (E1): $x-t-$diagram for $\varrho$
    \label{fig:E1_xt}
    }
\end{figure}

The density distribution results from an optimal excitation (see figure~\ref{fig:E1_xt} (left))  in pipe \RM{1}, which is transported through pipe \RM{2}, see figure~\ref{fig:E1_xt} (right).

Accordingly, the previously presented discretization is suitable to model the fluid mechanical processes in a network.
Besides, the adjoint-based on this discretization proves to be able to carry out a typical optimization. 
For the sake of brevity, a validation of the adjoint gradient is omitted.


\subsection*{(E1.2) Three pipes in a row - Steep gradient}

With the same three pipe setup, a steep gradient is tested as well. The initial condition is a jump in the velocity from five to zero in pipe \RM{1}. With that, a wave travels to the left out of the network and a shock to the right. With the filtering, the shock moves smoothly from one pipe to the next.

\begin{figure}
    \includegraphics[width = 0.49\textwidth]{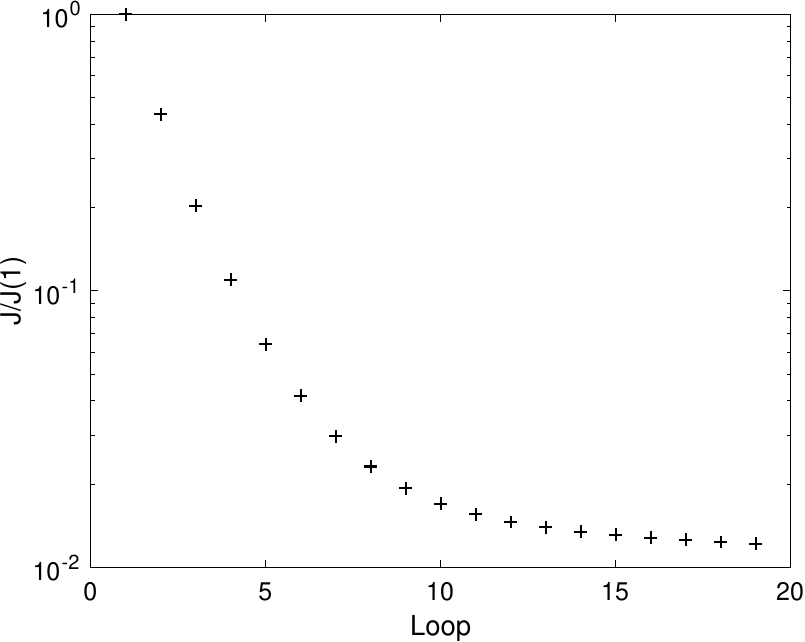}
    \includegraphics[width = 0.49\textwidth]{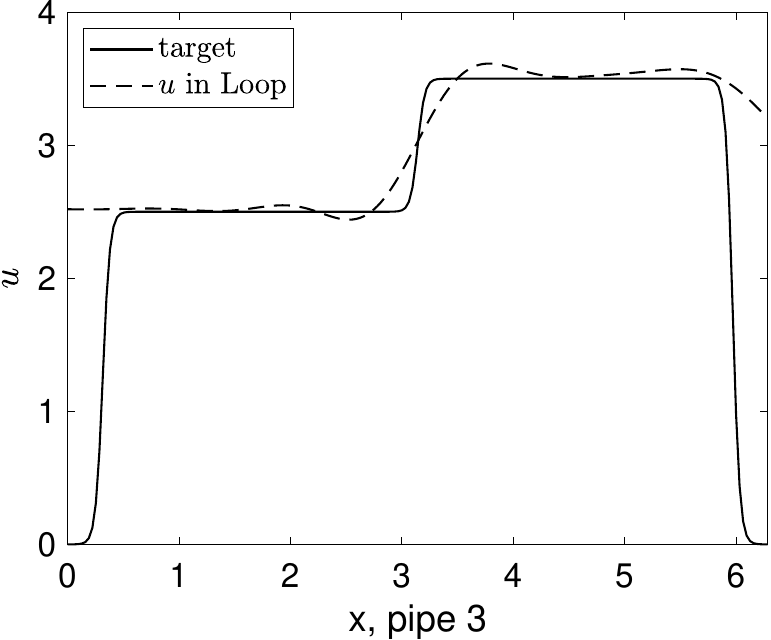}
    \caption{(Left) Progress of the objective function of (E1.2) normalized with respect to initial loop $J/J_1$. (Right) Velocity profile $u$ at the penultimate time step after 19 iterative loops.
    \label{fig:E1_2_J}
    }
\end{figure}

The defined goal in the penultimate time step is shown in figure~\ref{fig:E1_2_J} (right). It is a jump in the velocity in pipe~\RM{3} from 2.5 to 3.5.

The convergence criterion
\begin{equation}
    \frac{J_{n-1}-J_n}{J_1} < 2\cdot 10^{-4}  
\end{equation}
is reached after 19 loops. The objective function normalized by $J_1$ is reduced by nearly two orders of magnitude (see figure~\ref{fig:E1_2_J} (left)).

\begin{figure}
    \includegraphics[width = 0.49\textwidth]{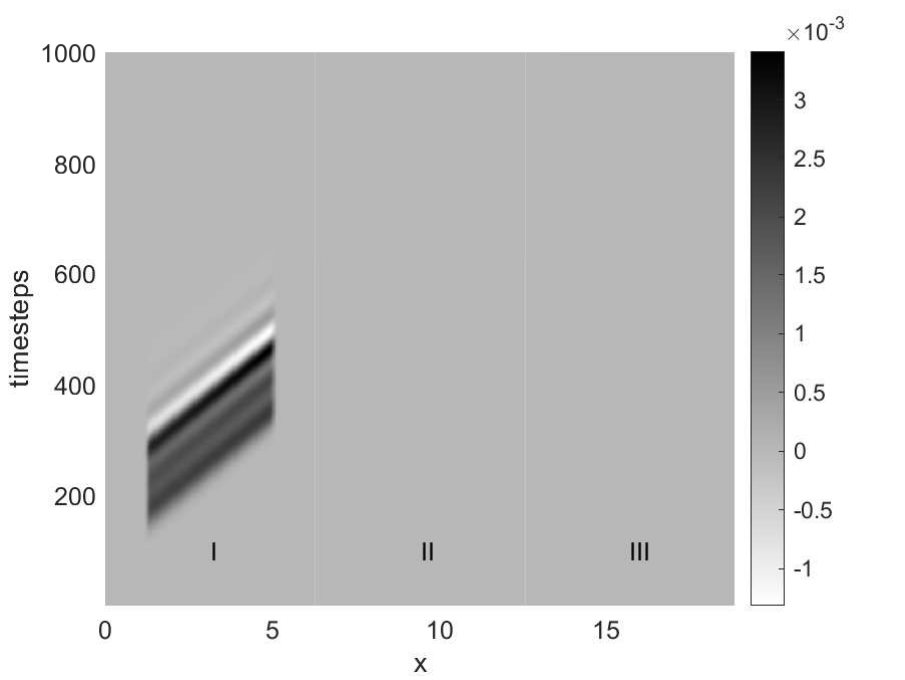}
    \includegraphics[width = 0.49\textwidth]{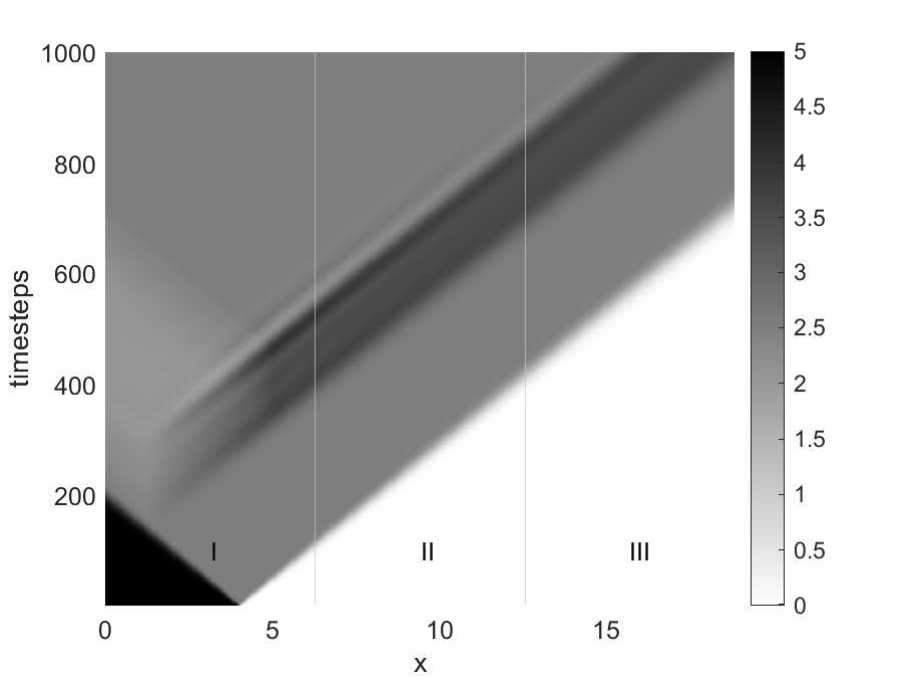}
    \caption{(Left) (E1) $x-t-$diagram for the excitation $r$. (Right) (E1): $x-t-$diagram for $u$
    \label{fig:E1_2_xt}
    }
\end{figure}

The velocity distribution results again from an optimal excitation (see figure~\ref{fig:E1_2_xt} (left))  in pipe \RM{1}. The gradient is flatter than the goal but the height of the jump is met.


\subsection{(E2) Diamond network} \label{sec:E2}

\begin{figure}
    \centerline{\includegraphics{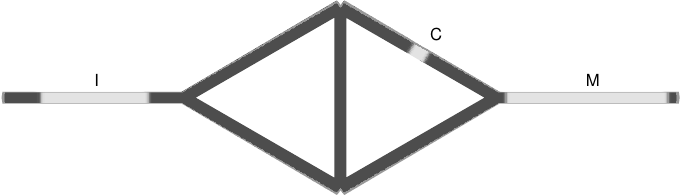}}
    \caption{Diamond network (E2) with M: measurement, C: consumer and I: influence}
    \label{fig:E2_setup}
\end{figure}

The second example considers a diamond-shaped network to mimic a typical application scenario with a customer whose demand changes significantly over time, see figure \ref{fig:E2_setup}.

As in the previous example the objective function is defined in terms of the density \eqref{eq_examples_objective}.
The reference density $\varrho_{target}$ is given by a constant density in pipe \RM{7}.

The overall goal is to ensure a uniform mass flow at the end of the network (M:measurement), taking into account a consumer (C:consumer) through optimal control of the mass flow at the beginning of the network (I:influence).
The control represented by additional terms $r_{\varrho}(x,t)$ and $r_m(x,t)$ on the right-hand-side of the momentum equation mimicking a controllable gas supply.
The consumer is modeled in the same way by a predefined disturbance $f_{disturb}$ in the momentum equation which is prescribed by
\begin{align}
    f_{disturb} = \frac{1}{2} \left( \tanh \left( \frac{1}{2} \left( t-0.1 \tau \right) \right) - \tanh \left( \frac{1}{2} \left( t-0.3\tau \right) \right) \right) \no\\+ \frac{1}{2} \left( \tanh \left( \frac{1}{2} \left( t-0.5\tau \right) \right) - \tanh \left( \frac{1}{2} \left( t-0.55\tau \right) \right) \right),
    \label{eq:disturbanceE2}
\end{align}
with $\tau$ being the total simulation time. This means the force is applied on the system from $10$ to $30~\%$ and from $50$ to $55~\%$ of the computational time. A pre-factor to equation~\eqref{eq:disturbanceE2} is chosen to see a significant change in the behavior of the system.

The governing parameters of the setup are stated in table~\ref{tab:parameter}.

\begin{figure}
    \includegraphics[width = 0.48\textwidth]{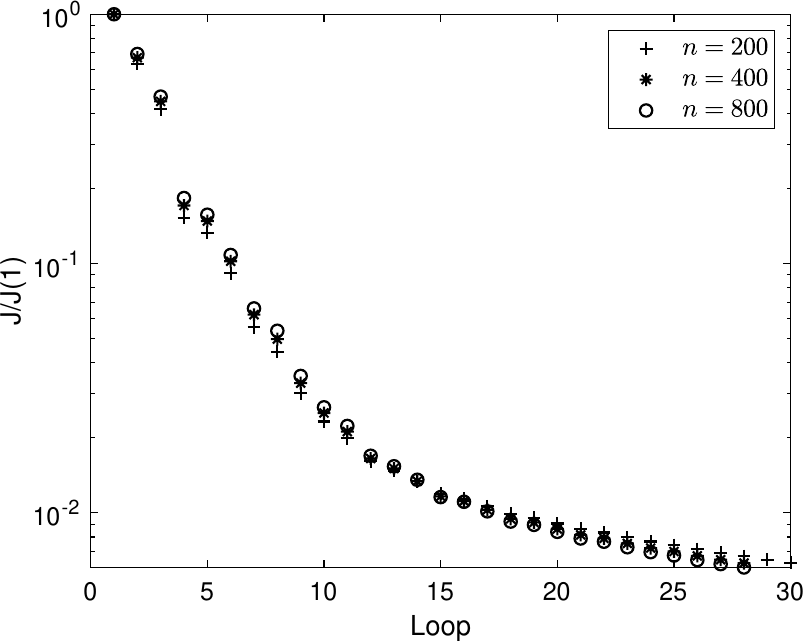}
    \includegraphics[width = 0.49\textwidth]{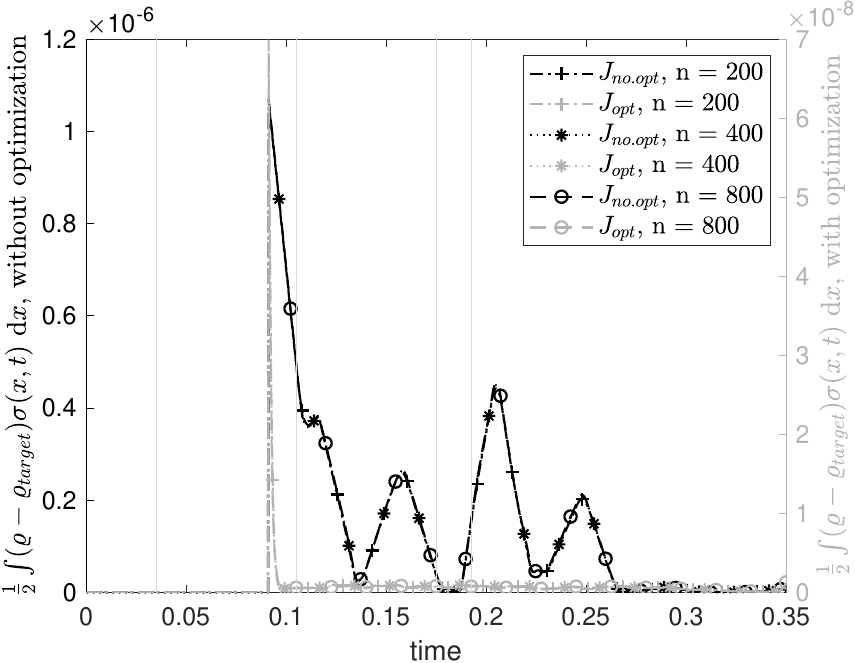}
    \caption{(Left) Progress of the objective function of (E2) normalized with respect to initial loop $J/J_1$. (Right) Spatial part of the integral of the objective over time after 30  and 28 iterative loops respectively, both for three different grid resolutions}.
    \label{fig:E2_J}
\end{figure}

Figure~\ref{fig:E2_J} (left) shows the decreasing objective function for the different grid resolutions. The convergence criterion of
\begin{equation}
    \frac{J_{n-1}-J_n}{J_1} < 2\cdot 10^{-4}
\end{equation}
is reached after 30 iterations for $200$ grid points per pipe and after 28 iterations for the higher resolutions of $n=400$ and $n=800$. The slope of the objective function shows that a further decrease is possible but not necessary as the optimization goal in pipe \RM{7} is already reached in typically sufficient precision.

\begin{figure}
    \centerline{\includegraphics[width = 0.49\textwidth]{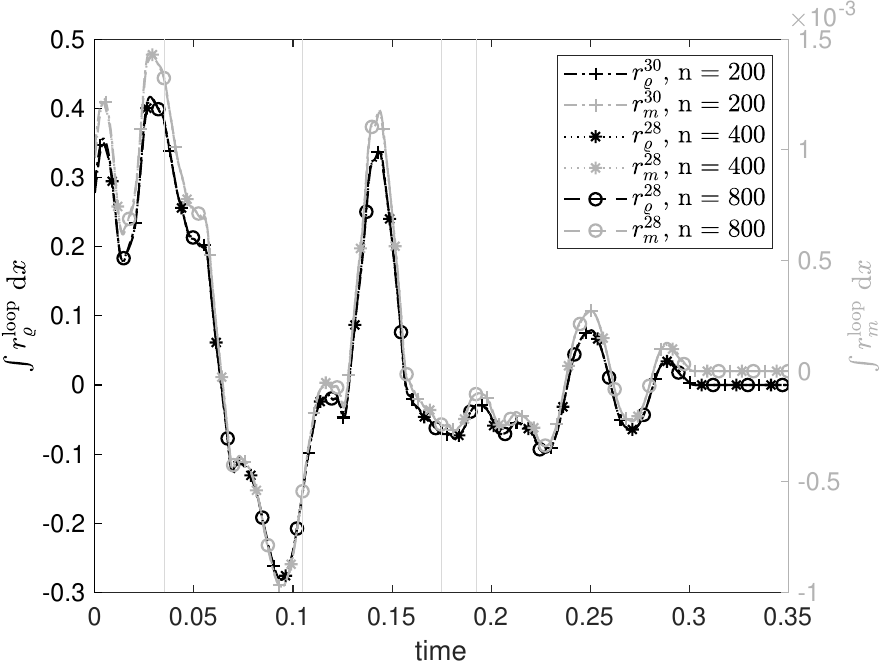}}
    \caption{Integrals of $r_{\varrho}$ and $r_m$ in pipe \RM{1} after 30 and 28 iterations depending on the grid resolution.
    \label{fig:E2_influence}
    }
\end{figure}

For the first $1300$ time steps the objective function $J$ is zero due to the choice of $\sigma$ as the adjoint does not influence the result in that time interval due to the distance between pipe~\RM{1} and the measurement in pipe~\RM{7}. The objective function is reduced by two orders of magnitude. The difference between the density $\varrho$ and $\varrho_{target}$ can be reduced significantly which leads to the spatial part of the integral of the objective function shown in figure~\ref{fig:E2_J} (right). Please note the different orders of magnitude between the left (without optimization) and the right (with optimization) y-axis. The differences between the different grid resolutions are negligible. The vertical lines indicate the two time intervals where the forcing is applied in pipe~\RM{5}. The value of the objective function is significantly decreased compared to the case without optimization. This is obtained by the $r_{\varrho}$ and $r_m$ in pipe \RM{1} shown in figure~\ref{fig:E2_influence}. Again, the differences between the different grids are negligible.


\begin{figure}
    \centerline{\includegraphics[width = \textwidth]{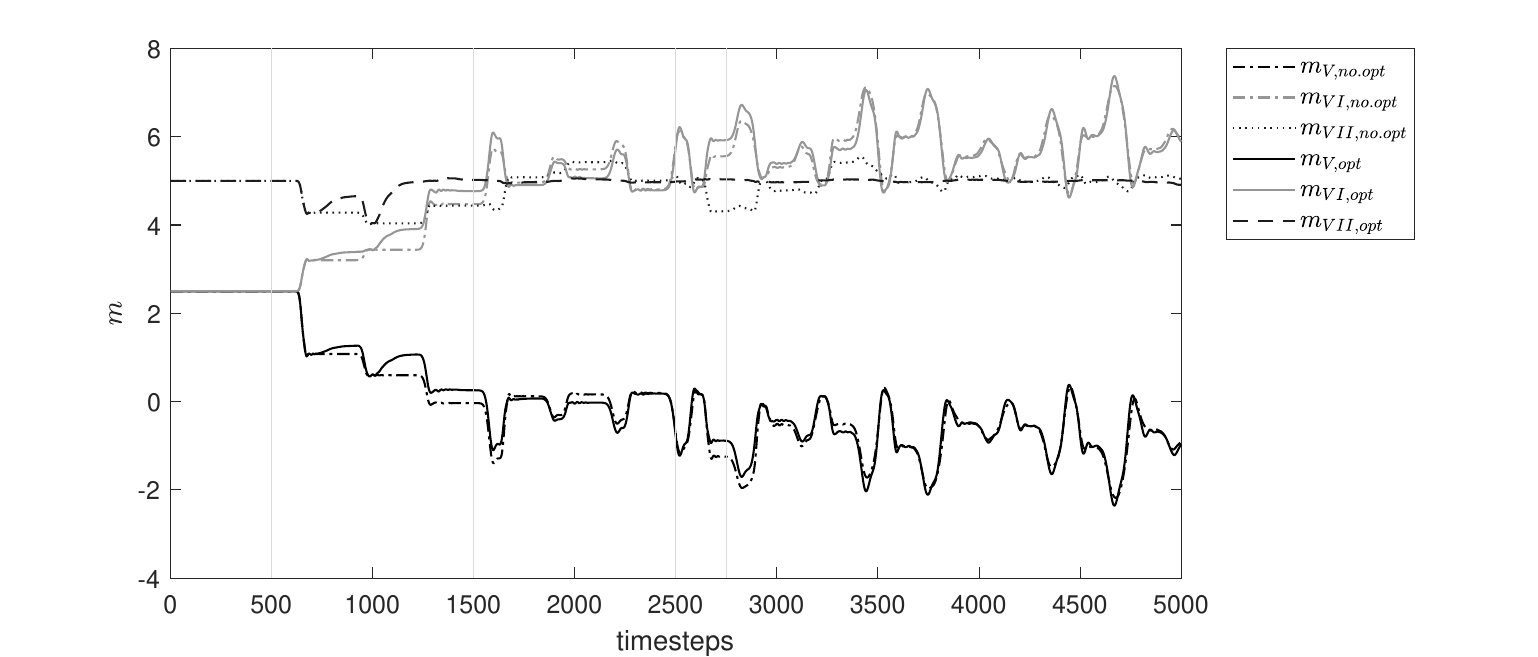}}
    \caption{Mass fluxes from pipe \RM{5} and \RM{6} to pipe \RM{7} with and without optimization ($n=200$).
    \label{fig:E2_massfluxes}
    }
\end{figure}

The mass fluxes from pipe \RM{5} and \RM{6} to pipe \RM{7} with and without optimization for the smallest grid resolution are shown in figure~\ref{fig:E2_massfluxes}. Again, the vertical lines mark the time intervals with disturbance. It can be observed that $m$ decreases in pipe~\RM{5} while there is a smaller increase in pipe~\RM{6} leading to a decrease in the total mass flux in pipe~\RM{7}. This first drop in $m$ cannot be compensated by the adjoint because the information needs some time to reach the last pipe. But after $1300$ time steps a nearly constant state is reached comparing the dashed and the dotted lines for $m_{VII}$ with and without optimization.

Also for this application-oriented configuration, the presented finite-difference discretization provides appropriate results.
The adjoint again provided suitable gradient information to substantially reduce the target function and to make the flow in the measuring area considerably more uniform.


\section{Conclusion}

We introduced a new finite differences approach to simulate gas networks. With the summation by part characteristic, it is possible to define the fluxes between the pipes adequately. The coupling conditions are the conservation of mass and equal pressure in all pipes at one node. The structure of the approach makes it possible to derive a simple formulation for the adjoint of the network, which allows for optimization tasks. 

In two sample networks, three connected pipes in a row and a diamond network,  it is shown that the developed framework is applicable, and the adjoint approach optimizes the system state towards the target defined in the objective function.

\section*{Acknowledgments}
The author acknowledges financial support by the Deutsche Forschungsgemeinschaft (DFG) within the collaborative Research Center (SFB) 1029 (200291049).

\bibliography{local}{}
\bibliographystyle{plain}


\appendix

\section{Complement to the Norm} \label{appendix:norm}

This section is in addition to the norm in section \ref{sec:numerical_implementation}.
For a symmetric and positive definite matrix $W$ the norm is defined as
\begin{equation}
 <u,v>_W= u^T W v .  
\end{equation}
With the used relations at the right side it follows
\begin{align}
    <u,D_x v>_W &= u^T W D_x v   &|& WD_x = S \\ 
    &= \left(\left(u^T S v \right)^T\right)^T    \\
    &= \left(v^T S^T u \right)^T               &|& S^T = -S + B\\
    &= \left(-v^T S u + v^T B u \right)^T \\
    &= -\left(v^T W D_x u \right)^T + \left(v^T B u \right)^T \\
    &= - \left(D_x u \right)^T W^T v + u^T B^T v  &|& W^T = W, B^T = B, B_{i,j} = -\delta_{i,1}\delta_{j,1} + \delta_{i,N}\delta_{j,N} \\
    &= -<D_x u, v>_W - u_1 v_1 + u_N v_N
\end{align}

This is the discrete analogous to partial integration.

\section{Non-reflecting boundary conditions}

The direct and adjoint non-reflecting boundary conditions are based on a characteristic decomposition of the flow field and a corresponding projection \cite{Thompson1987,PoinsotLele1992}.
While outgoing waves are defined within the computational domain, incoming waves are chosen to introduce no additional disturbances.

\paragraph{Direct}
For the derivation of the direct non-reflecting boundary conditions the governing equations \eqref{eq_euler_conti}-\eqref{eq_euler_momentum} are written in quasi-linear form 
\begin{equation}
    \begin{pmatrix} 1 & 0 \\ 0 & 1 \end{pmatrix} \partial_t q + \begin{pmatrix} 0 & 1 \\ 0 & 1 \end{pmatrix} \partial_x q = \partial_t q + \underbrace{\begin{pmatrix}0 & 1\\c^{2} - \frac{m^{2}}{\varrho^{2}} & \frac{2 m}{\varrho}\end{pmatrix}}_{Z} \partial_x q = 0
\end{equation}
with $q = [\varrho,m]$ and $m = \varrho u$.
The eigenvalues of the resulting operator $Z$ are $\lambda_{\pm} = u \pm c$.
An eigendecomposition of $Z$ leads to
\begin{equation}
  \partial_t q + \underbrace{\begin{bmatrix} e_+ & e_- 	\end{bmatrix}}_{= T}  
		   \underbrace{\begin{bmatrix} u+c 	& 0 		\\ 0 	& u-c  \end{bmatrix}}_{= \Lambda}
		   \underbrace{\begin{bmatrix} e^+ 	\\ e^- 		\end{bmatrix}}_{= T^{-1}} \partial_{x}q = 0.
\end{equation}
The matrices $T$ and $T^{-1}$ contain the right $\underline e_j$ and the left $\underline e^j$ eigenvectors of $Z$, corresponding to the acoustic $(\pm)$ characteristics. 
Due to scaling bi-orthonormality $\underline e^j \cdot \underline e_i = \delta_{ij}$ holds.

The perturbation of the state $\Delta q = q - q_0$ at a computational boundary with respect to a reference state $q_0$ is projected onto the resulting characteristic eigenvector base.
Here, the reference state corresponds to the initial condition $q(t=0)$.
\begin{equation}
  \Delta q = \underbrace{e^\pm \Delta q}_{=r_\pm} e_\pm 
\end{equation}
According to the desired boundary condition the factors $r_\pm$ can be chosen depending on the direction of the corresponding waves.
For outgoing waves $r_\pm$ remain unchanged. 
To avoid incoming waves $r_\pm = 0$ is chosen.
The modified state 
\begin{equation}
  \tilde q = q_0 + \Delta q
\end{equation}
contains no more incoming disturbances in a linearised sense.
Thus, a non-reflecting condition is obtained.

\paragraph{Adjoint}
The adjoint equivalent for a non-reflecting boundary condition is a non-reflecting condition for the adjoint variables.
Otherwise reflections of the adjoint variables would predict sensitivities on the objective function, which the direct system does not provide.

Thus, for adjoint non-reflecting boundary condition the same procedure can be applied.
The quasi-linear form of the adjoint equations is given by
\begin{equation}
    \partial_t q^* + \underbrace{\begin{pmatrix}0 & c^{2} - \frac{m^{2}}{\varrho^{2}}\\1 & \frac{2 m}{\varrho}\end{pmatrix}}_{Z} \partial_x q^* = 0
\end{equation}
resulting in the same eigenvalues $\lambda_{\pm} = u \pm c$ as in the direct governing equation which partially validates the adjoint equations.
Again a decomposition 
\begin{equation}
  \tilde q^* = q^*_0 + \Delta q^*
\end{equation}
and a suitable choice of the values $r^*\pm$ results in non-reflecting boundary conditions.
The reference value $q^*_0$ is chosen to the initial value $q^* = 0$.

\end{document}